\documentclass[aos,authoryear]{imsart}

\newcommand{\ignore}[1]{} 

\usepackage{amsthm,amsmath}

\usepackage[authoryear]{natbib}

\startlocaldefs

\newcommand{\mrm}[1]{\mathrm{#1}}
\newcommand{\ARL}{{\mrm{ARL2FA}}}
\newcommand{\ADD}{{\mrm{AD2D}}}

\newcommand{\mc}[1]{\mathcal{#1}} 
\newcommand{\Fc}{\mc{F}}
\newcommand{\Bc}{\mc{B}}

\newcommand{\mb}[1]{\boldsymbol{#1}} 
\newcommand{\Deb}{\mb{\Delta}}

\newcommand{\mbf}[1]{\mathbf{#1}} 
\newcommand{\Pb}{\mbf{P}}
\newcommand{\Eb}{\mbf{E}}

\newcommand{\de}{\delta}
\newcommand{\si}{\sigma}

\newcommand{\set}[1]{\left\{#1\right\}}

\newcommand{\brc}[1]{\left(#1\right)}
\newcommand{\brcs}[1]{\left[#1\right]}

\newcommand{\xra}{\xrightarrow}
\newcommand{\esup}{\operatornamewithlimits{ess\,sup}}

\newtheorem{theorem}{Theorem}

\newtheorem{corollary}{Corollary}

\theoremstyle{definition} 

\newtheorem{remark}{Remark}

\endlocaldefs

\begin{document}

\begin{frontmatter}

\title{On Optimality Properties of the Shiryaev-Roberts Procedure}
\runtitle{Optimality of the Shiryaev-Roberts Procedure}


\author{\fnms{Moshe} \snm{Pollak}\ead[label=e1]{msmp@mscc.huji.ac.il}\thanksref{t1}}
\thankstext{t1}{Moshe Pollak is Marcy Bogen Professor of Statistics at the Hebrew University of Jerusalem.
His work was supported in part by a grant from the Israel Science Foundation, by
the Marcy Bogen Chair of Statistics at the Hebrew University of Jerusalem, and by
the U.S.\ Army Research Office MURI grant W911NF-06-1-0094 at the University of
Southern California.}
\address{Hebrew University of Jerusalem\\
Department of Statistics\\
Mount Scopus\\
Jerusalem 91905, Israel\\
\printead{e1}}
 \and
\author{\fnms{Alexander G.} \snm{Tartakovsky}\corref{}\ead[label=e2]{tartakov@usc.edu}\thanksref{t2}\thanksref{t3}}
\thankstext{t2}{The work of Alexander Tartakovsky was supported in part by the U.S.\ Office of Naval
Research grant N00014-06-1-0110 and the U.S.\ Army Research Office MURI grant
W911NF-06-1-0094 at the University of Southern California.}
\thankstext{t3}{\textbf{Corresponding author.}}
\address{University of Southern California\\
Department of Mathematics\\
3620 S. Vermont Ave\\
Los Angeles, CA 90089-2532, USA\\
\printead{e2}}

\affiliation{The Hebrew University of Jerusalem and the University of Southern
California}

\runauthor{M. Pollak and A. G. Tartakovsky}

\begin{abstract}
We consider the simple changepoint problem setting, where observations are
independent, iid pre-change and iid post-change, with known pre- and post-change
distributions. The Shiryaev-Roberts detection procedure is known to be
asymptotically minimax in the sense of minimizing maximal expected detection delay
subject to a bound on the average run length to false alarm, as the latter goes to
infinity. Here we present other optimality properties of the Shiryaev-Roberts
procedure.

\end{abstract}



\end{frontmatter}

\section{Introduction} \label{s:Intro}

Changepoint problems deal with detecting a change in the state of a process, where
information one has about the state of affairs is in the form of observations. In
the sequential setting, observations are obtained sequentially and, as long as
their behavior is consistent with the initial (or target) state, one is content to
let the process continue. If the state changes, then one is interested in detecting
that a change is in effect, usually as soon as possible after its occurrence.

Any detection policy may give rise to false alarms. Intuitively, the desire to
detect a change quickly causes one to be (relatively) trigger-happy, which will
bring about many false alarms if there is no change. On the other hand, attempting
to avoid false alarms too strenuously will lead to a long delay between the time of
occurrence of a real change and its detection. Common operating characteristics of
a sequential detection policy are ARL2FA = the Average Run Length (the expected
number of observations) to False Alarm (assuming that there is no change) and the
AD2D = Average Delay to Detection (the expected delay between a real change and its
detection). The gist of the changepoint problem is to produce a detection policy
that (at least approximately) minimizes the AD2D subject to a bound on the ARL2FA.
The constitution of a good policy depends very much on what is known about the
stochastic behavior of the observations, both pre- and post-change.

Let $X_1, X_2, \dots$ denote the series of observations, and let $\nu$ be the
serial number of the first post-change observation. Let $\Pb_k$ and $\Eb_k$ denote
probability and expectation when $\nu=k$, and let $\Pb_\infty$ and $\Eb_\infty$
denote the same when $\nu=\infty$ (i.e., there never is a change). A sequential
change detection procedure is identified with a stopping time $N$ on
$X_1,X_2,\dots$, i.e., $\{N \le n\} \in \Fc_n$, where $\Fc_n=\si(X_1,\dots,X_n)$ is
the sigma-algebra generated by the first $n$ observations.

In this paper, we consider the simplest setting of the problem, where the
observations are independent, each having density $f_0$ pre-change and density
$f_1$ post-change, where both $f_0$ and $f_1$ are known, and only the value of
$\nu$, the point of change, is unknown. (In practice, often $f_0$ is known.
Realistically, $f_1$ is not known, but the simple setting yields a benchmark for
the best one can hope.) \ In this setting, \citet{MoustakidesAS86} proved that the
Cusum procedure \citep{Page54} is optimal in the sense of minimizing the
worst-worst case (essential supremum) expected detection delay
\[
\sup_{k \ge 1} \esup_{\Fc_k} \Eb_k[(N-k)^+|\Fc_k]
\]
over all stopping times $N$ for which
\begin{equation} \label{ARLconstraint}
\ARL(N)=\Eb_\infty N \ge B,
\end{equation}
where $B>0$ is a value set before the surveillance begins. See also
\citet{LordenCPDAS71} and \citet{Ritov90}. For a continuous-time Brownian motion a
similar result has been established by \citet{Beibel96} and \citet{ShiryaevRMS96}.

\citet{PollakAS85} proved that the Shiryaev-Roberts procedure
\citep{ShiryaevTPA63,Roberts66} is asymptotically (as $B\to\infty$) optimal in the
sense of minimizing the supremum AD2D
\[
\sup_{k \ge 1} \Eb_k(N-k|N \ge k)
\]
over all stopping times $N$ that satisfy \eqref{ARLconstraint}.

Here we prove other (exact) optimality properties of the Shiryaev-Roberts detection
procedure. To be specific, in Section~\ref{s:ExactOpt}, we prove that the
Shiryaev-Roberts procedure is (exactly) optimal in the sense of minimizing the
``integral  $\ADD"=\sum_{k=1}^\infty\Eb_k(N-k)^+$ for every $B>0$ in the class of
procedures with the ARL2FA constraint \eqref{ARLconstraint}. In
Section~\ref{s:Stationary}, we consider the setting where a change occurs in a
distant future (i.e., $\nu$ is large) and the detection of a change is preceded by
a large number of false detections. We prove that the Shiryaev-Roberts procedure is
the best one can do in terms of minimizing the expected detection delay
asymptotically when $\nu \to \infty$ in the class \eqref{ARLconstraint}, for every
$B>0$.

Both problem settings have been previously considered for a continuous-time
Brownian motion model. See \citet{FeinbergShiryaev07,ShiryaevTPA63} and
Remarks~\ref{Rem1} and \ref{Rem3} below.

\section{Minimizing integral AD2D} \label{s:ExactOpt}

Using the notation of the previous section, the Shiryaev-Roberts procedure calls
for stopping and raising an alarm at
\begin{equation} \label{SRst}
N_{A_B}= \min \set{n \ge 1: R_n \ge A_B},
\end{equation}
where
\begin{equation} \label{SRstat}
R_n = \sum_{k=1}^n \frac{p(X_1,\dots,X_n|\nu=k)}{p(X_1,\dots,X_n|\nu=\infty)}=
\sum_{k=1}^n \prod_{i=k}^n \frac{f_1(X_i)}{f_0(X_i)}
\end{equation}
and $A_B$ is such that $\Eb_\infty N_{A_B}=B$.

Below in Theorem~\ref{Th1} we prove that the Shiryaev-Roberts procedure is exactly
optimal in the sense of minimizing the integral $\ADD=\sum_{k=1}^\infty
\Eb_k(N-k)^+$ in the class of detection procedures $\Deb_B=\{N: \ARL(N) \ge B\}$ in
which the mean time to false alarm is not less than the given positive number $B$.
We begin with a sketch of the argument why one may expect this to be true.

To this end, we first need to consider the following Bayesian problem, denoted by
$\Bc(\rho,c)$. Suppose $\nu$ is random and has a geometric prior distribution
\[
\Pb(\nu=k)=\rho(1-\rho)^{k-1}, \quad k \ge 1,
\]
and the losses associated with stopping at time $N$ are 1 if $N< \nu$ and $c \cdot
(N-\nu)$ if $N \ge \nu$, where $0<\rho <1$ and $c>0$ are fixed constants. Write
$\Pb^\rho(\bullet)=\sum_{k=1}^\infty \rho(1-\rho)^{k-1}\Pb_k(\bullet)$ for the
``average" probability and $\Eb^\rho$ for the corresponding expectation.

Solution of $\Bc(\rho,c)$ requires minimization of the expected loss
\begin{equation} \label{ExpLoss}
\varphi_{c,\rho}(N)= \Pb^\rho(N<\nu)+ c \Eb^\rho(N-\nu)^+ ,
\end{equation}
and the Bayes rule for this problem is given by the Shiryaev procedure (cf.
\citealp{ShiryaevTPA63,Shiryaevbook78}), which is the stopping time
\begin{equation} \label{Shirst}
T_{\rho,c}= \min \set{n \ge 1: \Pb^\rho(\nu \le n|\Fc_n) \ge \de_{\rho,c}},
\end{equation}
where $0<\de_{\rho,c}<1$ is an appropriate threshold.

Obviously, the $\Bc(\rho,c)$ problem is equivalent to maximizing
\[
\frac{1}{\rho}[1-\varphi_{c,\rho}(N)]= \frac{\Pb^{\rho}(N\ge \nu)}{\rho} - c
\frac{\Eb^{\rho}(N-\nu)^+}{\rho}.
\]
In the proof of Theorem~\ref{Th1} below, we show that, for any stopping time $N$,
\[
\frac{\Pb^{\rho}(N\ge \nu)}{\rho} \xra[\rho \to 0]{} \Eb_\infty N, \quad
\frac{\Eb^{\rho}(N-\nu)^+}{\rho} \xra[\rho \to 0]{} \sum_{k=1}^\infty \Eb_k(N-k)^+
.
\]
Hence
\[
\frac{1}{\rho}[1-\varphi_{c,\rho}(N)] \xra[\rho \to 0]{}\Eb_\infty N -c
\sum_{k=1}^\infty \Eb_k(N-k)^+,
\]
which should be maximized in the class $\Deb_B$.

We also show that the Shiryaev procedure $T_{\rho,c}$ converges to the
Shiryaev-Roberts procedure $N_{A_B}$ as $\rho \to 0$. Therefore, it stands to
reason that the integral $\ADD=\sum_{k=1}^\infty \Eb_k(N-k)^+$ is minimized subject
to $\Eb_\infty N\ge B$.

Formal details are given in the following theorem and its proof.

\begin{theorem} \label{Th1}
Let $A_B$ be chosen so that $\ARL(N_{A_B})=B$. Then the Shiryaev-Roberts procedure
defined by \eqref{SRst} and \eqref{SRstat} minimizes
\begin{equation} \label{IntADD}
\sum_{k=1}^\infty \Eb_k(N-k)^+
\end{equation}
over all stopping times $N$ that satisfy $\Eb_\infty N \ge B$, i.e.,
\[
\inf_{N\in \Deb_B} \sum_{k=1}^\infty \Eb_k(N-k)^+ = \sum_{k=1}^\infty
\Eb_k(N_{A_B}-k)^+ \quad \text{for every $B >0$},
\]
where $\Deb_B=\{N: \ARL(N) \ge B\}$.
\end{theorem}

\begin{proof}
Consider the Bayesian problem $\Bc(\rho,c)$ with Geometric($\rho$) prior
distribution and the average loss \eqref{ExpLoss}.
\citet{ShiryaevTPA63,Shiryaevbook78} proved that the expected loss  \eqref{ExpLoss}
for the problem $\Bc(\rho,c)$ is minimized by the stopping time \eqref{Shirst}.
Applying Bayes' formula, it is easy to see that
\[
\Pb^\rho(\nu \le n|\Fc_n) = \frac{R_{\rho,n}}{R_{\rho,n} + 1/\rho},
\]
where
\[
R_{\rho,n} = \sum_{k=1}^n \prod_{i=k}^n \brc{\frac{1}{1-\rho}
\frac{f_1(X_i)}{f_0(X_i)}}.
\]
Hence, the Shiryaev rule can be written in the equivalent form
\begin{equation} \label{Shirst1}
T_{\rho,c}= \min \set{n \ge 1: R_{\rho,n} \ge A_{\rho,c}},
\end{equation}
where $A_{\rho,c}=(1/\rho)[\de_{\rho,c}/(1-\de_{\rho,c})]$.

Note first that $R_{\rho,n} \xra[\rho\to 0]{} R_n$.

By Theorem 1 of \citet{PollakAS85}, there exist a constant $0<c^*<\infty$ and a
sequence $\{\rho_i,c_i\}_{i=1}^\infty$ with $\rho_i\xra[i\to\infty]{}0$,
$c_i\xra[i\to\infty]{} c^*$ such that $N_{A_B}$ is the limit of the Bayes rules
$T_{\rho_i,c_i}$ as $i\to \infty$. Furthermore,
\begin{equation} \label{LimSup}
\limsup_{p \to 0, c\to c^*}
\frac{1-\varphi_{c,\rho}(T_{\rho,c})}{1-\varphi_{c,\rho}(N_{A_B})} =1,
\end{equation}
where $\varphi_{c,\rho}(N)$ is the expected loss associated with using the stopping
time $N$ for $\Bc(\rho,c)$.

Now, for any stopping time $N$,
\begin{equation*} 
\begin{split}
\frac{1}{\rho}[1-\varphi_{c,\rho}(N)] & = \frac{1}{\rho} \brcs{(1-\Pb^\rho(N<\nu))
- c \Eb^\rho (N-\nu)^+}
\\
& = \frac{\Pb^\rho(N\ge \nu)}{\rho} \brcs{1- c \Eb^\rho (N-\nu|N\ge \nu)}.
\end{split}
\end{equation*}
Since
\begin{equation*} 
\begin{split}
\frac{\Pb^\rho(N\ge \nu)}{\rho} & = \frac{1}{\rho} \sum_{k=1}^\infty \Pb_k (N \ge
k) \rho(1-\rho)^{k-1}
\\
& = \sum_{k=1}^\infty \Pb_\infty (N \ge k) (1-\rho)^{k-1}
\\
& \xra[\rho \to 0]{} \sum_{k=1}^\infty \Pb_\infty (N \ge k) =\Eb_\infty N
\end{split}
\end{equation*}
and
\begin{align*}
\frac{\Pb^\rho(N\ge \nu)\Eb^\rho(N-\nu|N\ge \nu)}{\rho} & = \frac{\Eb^\rho(N-\nu;
N\ge \nu)}{\rho}
\\
& = \frac{1}{\rho}\sum_{k=1}^\infty \Eb_k(N-k; N\ge k) \rho(1-\rho)^{k-1}
\\
& = \sum_{k=1}^\infty \Eb_k(N-k; N\ge k) (1-\rho)^{k-1}
\\
& \xra[\rho \to 0]{} \sum_{k=1}^\infty \Eb_k(N-k; N\ge k)= \sum_{k=1}^\infty
\Eb_k(N-k)^+ ,
\end{align*}
it follows that for any stopping time $N$ that has finite ARL2FA
\[
\frac{1}{\rho}[1-\varphi_{c,\rho}(N)] \xra[\rho \to 0]{} \Eb_\infty N - c
\sum_{k=1}^\infty \Eb_k(N-k)^+,
\]
which together with \eqref{LimSup} establishes that the Shiryaev-Roberts procedure
minimizes \eqref{IntADD} over all stopping times that satisfy $\Eb_\infty N=B$.
Note that if $B_1 >B$, then $N_{A_{B_1}}$ is stochastically larger than $N_{A_B}$,
i.e., all expectations in \eqref{IntADD} become larger. This implies that the
Shiryaev-Roberts procedure minimizes \eqref{IntADD} in the class $\Deb_B$. This
completes the proof of the theorem.
\end{proof}

\begin{corollary} \label{Cor1}

The Shiryaev-Roberts procedure defined by \eqref{SRst} and \eqref{SRstat} minimizes
\begin{equation} \label{ADDCondInt}
\frac{\sum_{k=1}^\infty \Eb_k(N-k|N\ge k)\Pb_\infty(N\ge k)}{\sum_{j=1}^\infty
\Pb_\infty(N\ge j)}
\end{equation}
over all stopping times $N$ that satisfy $\Eb_\infty N = B$, i.e.,
\[
\inf_{\{N: \Eb_\infty N =B\}} \sum_{k=1}^\infty w_k(N)\Eb_k(N-k|N\ge k) =
\sum_{k=1}^\infty w_k(N_{A_B})\Eb_k(N_{A_B}-k|N_{A_B}\ge k),
\]
where
\[
w_k(N) = \frac{\Pb_\infty(N\ge k)}{\sum_{j=1}^\infty \Pb_\infty(N\ge j)}
\]
and the threshold $A_B$ is selected so that $\Eb_\infty N_{A_B}=B$.
\end{corollary}

\begin{proof}

Obviously,  $\sum_{j=1}^\infty \Pb_\infty(N\ge j)=\Eb_\infty N =B$, so the
denominator in \eqref{ADDCondInt} is constant over all stopping times under
consideration. As for the numerator,
\begin{equation} \label{ADDk}
\begin{split}
\Eb_k(N-k|N\ge k)\Pb_\infty(N\ge k) & = \Eb_k(N-k|N\ge k)\Pb_k(N\ge k)
\\
& =\Eb_k(N-k; N\ge k) = \Eb_k(N-k)^+.
\end{split}
\end{equation}
Application of Theorem \ref{Th1} concludes the proof.
\end{proof}

\begin{remark} \label{Rem1}
Recently, \citet{FeinbergShiryaev07} established a result similar to
Theorem~\ref{Th1} for Brownian motion where an abrupt change occurs in the drift,
in which case the integral AD2D is $\int_0^\infty \Eb_\nu (N - \nu)^+ d\nu$. They
refer to this as ``A Generalized Bayesian Setting."
\end{remark}

While Theorem~\ref{Th1} and Corollary~\ref{Cor1} are of interest in their own
right, they are useful for proving other interesting optimality results, as it will
become apparent in the next section.

\section{Optimality for a change appearing after many re-runs} \label{s:Stationary}

Consider a context in which it is of utmost importance to detect a real change as
quickly as possible after its occurrence, even at the price of raising many false
alarms (using a repeated application of the same stopping rule) before the change
occurs. This essentially means that the changepoint $\nu$ is very large compared to
the constant $B$ which, in this case, defines the mean time between consecutive
false alarms.

To be more specific, let $N_{A_B}^{(1)}, N_{A_B}^{(2)}, \dots$ be sequential
independent repetitions of $N_{A_B}$ defined in \eqref{SRst}, i.e.,
\begin{equation} \label{Ni}
N_{A_B}^{(i)} = \min \set{n > \sum_{j=1}^{i-1} N_{A_B}^{(j)}: R_n^{(i)} \ge A_B} -
\sum_{j=1}^{i-1} N_{A_B}^{(j)} ,
\end{equation}
where $N_{A_B}^{(0)}=0$ and
\begin{equation} \label{SRi}
R_n^{(i)} = \sum_{k=N_{A_B}^{(i-1)}+1}^n \prod_{i=k}^n \frac{f_1(X_i)}{f_0(X_i)},
\quad \sum_{j=1}^{i-1} N_{A_B}^{(j)} < n \le \sum_{j=1}^{i} N_{A_B}^{(j)} .
\end{equation}
Therefore, $R_n^{(i)}$, $n >  \sum_{j=1}^{i-1} N_{A_B}^{(j)}$  is nothing but the
Shiryaev-Roberts statistic that is renewed from scratch after the $(i-1)$st false
alarm (under $\Pb_\infty$) and is applied to the segment of data
\[
X_{\sum_{j=1}^{i-1}N_{A_B}^{(j)}+1}, X_{\sum_{j=1}^{i-1}N_{A_B}^{(j)}+2}, \dots .
\]
Note that $\Eb_\infty N_{A_B}^{(i)} =B$ for $i \ge 1$.

Let, for $j \ge 1$,
\begin{equation} \label{tauj}
Q_j = N_{A_B}^{(1)}+ N_{A_B}^{(2)} + \cdots + N_{A_B}^{(j)}
\end{equation}
be the time of the $j$-th alarm, and let $J_\nu =\min\{j \ge 1: Q_j\ge \nu\}$,
i.e., $Q_{J_\nu}$ is the time of detection of a true change that occurs at $\nu$
after $J_\nu-1$ false alarms have been raised.

Our next theorem states that the Shiryaev-Roberts procedure defined by $Q_{J_\nu}$
is asymptotically (as $\nu \to \infty$) optimal with respect to the expected delay
$\Eb_\nu(Q_{J_\nu}-\nu)$ in the class of detection procedures for which the mean
time between false alarms is not less than $B$. Note that this result is not
asymptotic with respect to the ARL2FA. In fact, it holds for every positive $B$.

\begin{theorem} \label{Th2}
Let $\nu$ be the time of the change. Let  $N_{A_B}^{(1)}, N_{A_B}^{(2)}, \dots$ be
sequential independent repetitions of $N_{A_B}$ as defined in \eqref{Ni} and let
$Q_1,Q_2, \dots$ be as in \eqref{tauj}. Let $J_\nu =\min\{j: Q_j\ge \nu\}$.

(i) $\lim_{\nu \to \infty} \Eb_\nu(Q_{J_\nu}-\nu)$ exists.

(ii) Suppose a detection procedure $N$ with $\ARL(N)=B$ is applied repeatedly. Let
$N_1, N_2, \dots$ be sequential repetitions of $N$, let $W_j = \sum_{i=1}^j N_i$,
and let $K_\nu=\min\{j: W_j \ge \nu\}$. Then, for every $B>0$,
\begin{equation} \label{Optimality}
\lim_{\nu \to \infty} \Eb_\nu(Q_{J_\nu}-\nu) \le \lim_{\nu \to \infty}
\Eb_\nu(W_{K_\nu}-\nu).
\end{equation}

(iii) Inequality \eqref{Optimality} holds for all $N \in \Deb_B$, where
$\Deb_B=\{N: \Eb_\infty N \ge B\}$.
\end{theorem}

\begin{proof}
Proof of (i). By renewal theory, the distribution of $\nu-Q_{J_\nu-1}$ has a limit
\begin{equation} \label{Prob}
\lim_{\nu\to\infty} \Pb_\nu\brc{\nu-Q_{J_\nu-1}=k} = \frac{\Pb_\infty(N_{A_B}\ge
k)}{\sum_{j=1}^\infty \Pb_\infty(N_{A_B}\ge j)}
\end{equation}
(see, e.g., \citealp[page 356]{Feller-book66}).

Using \eqref{Prob} and letting $N_{A_B}$ be independent of $N_{A_B}^{(1)},
N_{A_B}^{(2)}, \dots$, we obtain
\begin{align*}
\Eb_\nu(Q_{J_\nu}-\nu) & = \Eb_\nu\Eb_\nu\brc{Q_{J_\nu}-\nu |Q_{J_\nu-1}}
\\
& = \sum_{k=1}^\nu \Eb_k\brc{N_{A_B}-k |\nu-Q_{J_\nu-1}=k, N_{A_B}\ge k} \Pb_\infty
\brc{\nu-Q_{J_\nu-1}=k}
\\
& = \sum_{k=1}^\nu \Eb_k\brc{N_{A_B}-k |N_{A_B}\ge k} \Pb_\infty
\brc{\nu-Q_{J_\nu-1}=k}
\\
& \xra[\nu \to \infty]{} \frac{\sum_{k=1}^\infty \Eb_k\brc{N_{A_B}-k |N_{A_B}\ge k}
\Pb_\infty \brc{ N_{A_B} \ge k}}{\sum_{j=1}^\infty\Pb_\infty \brc{ N_{A_B} \ge j}}
\\
& = \frac{\sum_{k=1}^\infty\Eb_k(N_{A_B}-k)^+}{\Eb_\infty N_{A_B}} =
\frac{\sum_{k=1}^\infty\Eb_k(N_{A_B}-k)^+}{B},
\end{align*}
which completes the proof of (i).

Proof of (ii). The same argument as in the proof of (i) yields
\[
\lim_{\nu \to \infty} \Eb_\nu(W_{K_\nu}-\nu) =
\frac{\sum_{k=1}^\infty\Eb_k(W-k)^+}{B}.
\]
Combining this with Corollary~\ref{Cor1} concludes the proof.

Proof of (iii). Write $\ADD(B)=\lim_{\nu \to \infty}\Eb_\nu(Q_{J_\nu}-\nu)$ for the
AD2D of the Shiryaev-Roberts procedure $N_{A_B}$. Note that $\ADD(B)$ tends to 0 as
$B \to 0$ and to $\infty$ as $B \to \infty$.  By virtue of (ii), it suffices to
show that $\ADD(B)$ is nondecreasing in $B$.

Note that $\ADD(B)$ is continuous in $B$. Therefore, if $\ADD(B)$ were not
nondecreasing in $B$, there would exist $0< B_1 <B_2 < \infty$ such that
$\ADD(B_1)=\ADD(B_2)$ and $\ADD(B) > \ADD(B_1)=\ADD(B_2)$ for all $B_1< B < B_2$.

Consider the following renewal-theoretic argument. Let $L_1, L_2, \dots$ and $M_1,
M_2, \dots$ be independent sequences of positive random variables having finite
means, each of them iid. Let $G^L$ be the asymptotic distribution of the residual
waiting time of the sequence $\{L_i\}$ (i.e., of the overshoot of the sequence
$\{\sum_{i=1}^j L_i\}_{j=1}^\infty$ over $t$, as $t \to \infty$) and let $G^M$ be
that of the sequence $\{M_i\}$. Let $G^T$ be the asymptotic distribution of the
residual waiting time for the sequence $\{T_i\}$ that is defined as follows:
$\Pb(T_i=L_i)=\Pb(T_i=M_i)=1/2$. By the usual renewal-theoretic apparatus, one can
show that

(a) $G^T = \tfrac{\Eb L}{\Eb L + \Eb M} G^N + \tfrac{\Eb M}{\Eb L + \Eb M} G^M$.

(b) Let $n_t = \min\{n: \sum_{i=1}^n T_i \ge t\}$. The asymptotic probability (as
$t\to \infty$) that $T_{n_t}$ is of the type $L,M$ is $\Eb L/(\Eb L+\Eb M)$, $\Eb
M/(\Eb L+\Eb M)$, respectively.

(c) Conditional on $T_{n_t}$ being the type $L,M$ the asymptotic (as $t\to \infty$)
distribution of the residual waiting time $\sum_{i=1}^{n_t} T_i -t$ is $G^L,G^M$,
respectively:
\[
\lim_{t\to \infty} \Pb\brc{\sum_{i=1}^{n_t} T_i -t \le x | T_{n_t} \; \text{is of
type}\; L, M} = G^L, G^M.
\]

Now, let $L=N_{A_{B_1}}$ and $M=N_{A_{B_2}}$. Note that the notation $n_t$ in terms
of $L$ and $M$ is the same as $J_t$ in terms of $N_{A_{B_1}}$ and $N_{A_{B_2}}$.
Recall that the procedure based on $T$ ``recycles" every time the Shiryaev-Roberts
statistic crosses the boundary ($A_{B_1}$ if the cycle has $T$ of type
$N_{A_{B_1}}$ and $A_{B_2}$ if the cycle has $T$ of type $N_{A_{B_2}}$). Let
$R^N_t$ be the value of the detection statistic at time $t$, where $N$ is a generic
stopping time that is applied repeatedly.

Let $R_n^{(i)}(N_{A_{B_1}})$ be equal to $R_n^{(i)}$ of \eqref{SRi} for $B=B_1$ and
let $R_n^{(i)}(N_{A_{B_2}})$ be the same for $B=B_2$. To emphasize the dependence
of $J_\nu$ and $Q_{j}$ on the stopping time $N$ being used, we will write $J_\nu^N$
and $Q^N(j)$. With this notation, for $j=1,2$,
\begin{align*}
& \Pb_\nu \brc{R^T_\nu \le x | Q^T(J_\nu^T-1)=\nu -k, T_{J_\nu^T} \; \text{is of
type} \; N_{A_{B_j}} }
\\
& = \Pb_\infty \brc{R_k^{(1)}(N_{A_{B_j}}) \le x | N_{A_{B_j}} \ge k }
\stackrel{\text{def}}{=}F_{j,k}(x).
\end{align*}
Therefore,
\begin{align*}
& \Pb_\nu \brc{R^T_\nu \le x}
\\
& = \frac{\Eb_\infty N_{A_{B_1}}}{\Eb_\infty N_{A_{B_1}} + \Eb_\infty N_{A_{B_2}}}
\sum_{k=1}^\infty F_{1,k}(x)
\\
& \quad \times \Pb_\infty \brc{Q^{N_{A_{B_1}}}(J_\nu^{N_{A_{B_1}}})=\nu-k |
T_{J_\nu^T} \; \text{is of type} \; N_{A_{B_1}}}
\\
& \quad + \frac{\Eb_\infty N_{A_{B_2}}}{\Eb_\infty N_{A_{B_1}} + \Eb_\infty
N_{A_{B_2}}} \sum_{k=1}^\infty F_{2,k}(x)
\\
& \quad \times \Pb_\infty\brc{Q^{N_{A_{B_2}}}(J_\nu^{N_{A_{B_2}}})=\nu-k |
T_{J_\nu^T} \; \text{is of type} \; N_{A_{B_2}}}
\end{align*}
and, by \eqref{Prob},
\begin{align*}
& \lim_{\nu\to\infty} \Pb_\nu \brc{R^T_\nu \le x}
\\
& = \frac{\Eb_\infty N_{A_{B_1}}}{\Eb_\infty N_{A_{B_1}} + \Eb_\infty N_{A_{B_2}}}
\sum_{k=1}^\infty F_{1,k}(x) \frac{\Pb_\infty \brc{N_{A_{B_1}} \ge k}}{\Eb_\infty
N_{A_{B_1}}}
\\
& \quad + \frac{\Eb_\infty N_{A_{B_2}}}{\Eb_\infty N_{A_{B_1}} + \Eb_\infty
N_{A_{B_2}}} \sum_{k=1}^\infty F_{2,k}(x) \frac{\Pb_\infty \brc{N_{A_{B_2}} \ge
k}}{\Eb_\infty N_{A_{B_2}}}.
\end{align*}

By abuse of notation, write $\ADD(N)$ for the limit (as $\nu\to \infty$) of the
average delay to detection when a stopping time $N$ is applied repeatedly. It now
follows that
\begin{align*}
\ADD(T)  & =\frac{\Eb_\infty N_{A_{B_1}}}{\Eb_\infty N_{A_{B_1}} + \Eb_\infty
N_{A_{B_2}}} \ADD(N_{A_{B_1}})
\\
&  \quad + \frac{\Eb_\infty N_{A_{B_1}}}{\Eb_\infty N_{A_{B_1}} + \Eb_\infty
N_{A_{B_2}}} \ADD(N_{A_{B_2}})
\\
& = \frac{B_1}{B_1 + B_2} \ADD(B_1) + \frac{B_2}{B_1 + B_2} \ADD(B_2)
\\
& = \ADD(B_1) = \ADD(B_2).
\end{align*}

Note that $\Eb_\infty T = \tfrac{1}{2} \Eb_\infty N_{A_{B_1}} +
\tfrac{1}{2}\Eb_\infty N_{A_{B_2}}= (B_1+B_2)/2$. By definition of $B_1$ and $B_2$,
it follows that
\[
\ADD(T) < \ADD(N_{A_{B}})=\ADD(B) \quad \text{for $B=\tfrac{1}{2}(B_1+B_2)$},
\]
which contradicts (ii) for $B=\tfrac{1}{2}(B_1+B_2)$.
\end{proof}

\begin{remark} \label{Rem2}
Theorem~\ref{Th2}(iii) implies that Corollary~\ref{Cor1} holds for all stopping
times $N\in \Deb_B$.
\end{remark}

\begin{remark} \label{Rem3}
\citet{ShiryaevTPA63} proved a result similar to Theorem~\ref{Th2}(ii) for Brownian
motion when a change occurs in the drift, and called this problem ``Quickest
Detection of a Disorder in a Stationary Regime."
\end{remark}

\begin{remark} \label{Rem4}
It is worth noting that Theorem~\ref{Th2} is important in a variety of surveillance
applications such as target detection and tracking, rapid detection of intrusions
in computer networks, and environmental monitoring, to name a few. In all of these
applications, it is of utmost importance to detect very rapidly changes that may
occur in a distant future, in which case the true detection of a real change may be
preceded by a long interval with frequent false alarms that are being filtered by a
separate mechanism or algorithm. For example, falsely initiated target tracks are
usually filtered by a track confirmation/deletion algorithm; false detections of
attacks in computer networks in anomaly-based Intrusion Detection Systems (IDS) may
be filtered by Signature-based IDS algorithms, etc. See, e.g.,
\citet{Tarbook91,TarVeerASMbook04,TartakovskyetalStamet06}. The practical
implication of Theorem~\ref{Th2} is that in these circumstances one has reason to
prefer the Shiryaev-Roberts procedure to other surveillance schemes.
\end{remark}

\end{document}